\documentclass{article}
\usepackage{amsmath,amssymb,amsthm}
\numberwithin{equation}{section}
\usepackage{verbatim}
\usepackage[toc,page]{appendix}
\usepackage{graphicx}
\usepackage{fourier}
\theoremstyle{plain} 
\newtheorem{thm}{Theorem}[section] 
\newtheorem*{theorem*}{Theorem}

\newtheorem{prop}[thm]{Proposition} 

\theoremstyle{definition}

\newtheorem*{defn*}{Definition}

\begin{document}
\title{A hydrodynamical homotopy co-momentum map  and a multisymplectic interpretation of higher order linking numbers}
\author{ Antonio Michele MITI${}^{\dagger, \dagger\dagger}$ and Mauro SPERA${}^{\dagger}$\\
\phantom{void}\\
${}^{\dagger}$ Dipartimento di Matematica e Fisica  ``Niccol\`o Tartaglia" \\ 
 Universit\`a Cattolica del Sacro Cuore\\
Via dei Musei 41, 25121 Brescia, Italia\\
\phantom{void}\\
${}^{\dagger\dagger}$ Departement Wiskunde, KU-Leuven\\ 
Celestijnenlaan 200B, B-3001 Leuven (Heverlee), Belgi\"e\\
\phantom{void}\\
{\small antoniomichele.miti@unicatt.it}, \,
{\small mauro.spera@unicatt.it} \\ }
\smallskip
\date{15th October 2019}
\maketitle
\begin{abstract}
In this article a homotopy co-momentum map (\`a la Callies-Fr\'egier-Rogers-Zambon) trangressing to the standard hydrodynamical co-momentum map of  Arnol'd, Marsden and Weinstein and others is constructed and then generalized to a special class of Riemannian manifolds. Also, a covariant phase space interpretation of the coadjoint orbits associated to the Euler evolution for perfect fluids and in particular of Brylinski's manifold of smooth oriented knots is discussed. 
As an application of the above homotopy co-momentum map, a reinterpretation of the (Massey) higher order linking numbers in terms of conserved quantities within the multisymplectic framework is provided and  knot theoretic analogues of first integrals in involution are determined.
\par
\end{abstract}
\smallskip



\section{Introduction}\par
In this paper we discuss some applications of  multisymplectic techniques in a hydrodynamical
context. 
The possibility of applying symplectic techniques therein ultimately comes from Arnol'd's pioneering work culminating in the geometrization of fluid mechanics (\cite{Arnold66, AM,Arn-Khe,MW83}).
In particular, in this connection we may mention the paper \cite{Rasetti-Regge75}, with its symplectic reinterpretation \cite{Pe-Spe89,Pe-Spe92,Pe-Spe00}, and the general portrait depicted in \cite{Bry}. 
Here we wish to apply some recently emerged concepts in multisymplectic geometry (mostly building on \cite{Calliesetal,RWZ,RW}) and 
construct an explicit {\it homotopy co-momentum map} (\cite{Calliesetal}) in a hydrodynamical setting, leading to
a multisymplectic interpretation of the so-called {\it higher order linking numbers}, viewed \`a la Massey (\cite{Pe-Spe02,Spe06,Hebda-Tsau12}). The costruction is generalized to cover connected compact oriented Riemannian manifolds having vanishing intermediate de Rham groups. Moreover, a {\it covariant phase space} intepretation of the multisymplectic setting is outlined.\par
We make clear from the outset that our constructions, together with the covariant phase space portrait, will not adhere to the standard multisymplectic approach to continuum mechanics set forth e.g. in \cite{Gotay&al, MPSW} but they will be based instead on the peculiar structure of an ideal fluid, whose configuration space is the ``Lie group" of  diffeomorphisms preserving a volume form, which will be directly taken as a multisymplectic form (\cite{CID}).\par
The layout of the paper is the following. First, in Section 2, we give an example of homotopy co-momentum map in fluid mechanics - in the sense of Callies-Fr\'egier-Rogers-Zambon (CFRZ) (\cite{Calliesetal}) 
- transgressing to Brylinski's symplectic structure on loop spaces and descending, in turn, to the manifold of smooth oriented knots, see \cite{Bry,BeSpe06} and below for precise definitions.  We briefly discuss  the (non) equivariance of the above construction with respect to the group of volume preserving diffeomorphisms of 3-space (see Section 2) and we outline a generalization thereof in a Riemannian framework, signalling  potential topological obstructions. Moreover,
covariant phase space aspects will be analyzed.
In Section 3 we prepare the ground for the forthcoming
applications by depicting a hydrodynamical multisymplectic portrait of basic knot theoretic objects, used,
in Section 4, to reinterpret the Massey higher order linking numbers in multisymplectic terms: the 1-forms appearing in the hierarchical Massey construction (viewed, in turn, differential geometrically \`a la Chen) provide an example of {\it first integrals in involution} in a multisymplectic framework. 
The last section is devoted to gathering together the conclusions and to pointing out possible directions for further research.
Appropriate background material is provided within the various sections in order to ease readability. \par
This paper is an improved version of part of the preprint \cite{Miti-Spera18}.\par

\section{Multisymplectic geometry and hydrodynamics of perfect fluids}\par
In the present section we freely use basic material on symplectic and multisymplectic geometry
tailored to our subsequent needs, prominently referring, for additional details, to \cite{Pe-Spe92,Spe06,Spera16} for the former and to \cite{RWZ,RW} for the latter. 
For general background on symplectic geometry and (co)momentum maps we quote, among others 
\cite{AM,GS84,Arn-Khe}.\par
\subsection{Tools in multisymplectic geometry}
All our objects will be smooth, unless differently specified. A (finite dimensional) {\it multisymplectic manifold} $(M, \omega)$ is a manifold (connected, for simplicity) equipped with a closed (n+1)-form $\omega$ 
(called {\it multisymplectic form} or {\it n-plectic form}) such that the map $\alpha$ below sending vector fields to n-forms (via contraction)
$$
{\mathfrak{X}}({M})  \ni \xi \mapsto  \alpha(\xi) := \iota_{\xi} \omega \in \Lambda^n(M)
$$
is injective (\cite{CID}). Dropping the last condition leads to the concept of {\it pre-n-plectic form}. The  $n=1$ case retrieves
(pre)symplectic manifolds.\par
In the multisymplectic context, the generalization of  the (co)momentum maps of the symplectic case leads to the more
refined concept of homotopy co-momentum map, to be presently succinctly reviewed.  One first introduces the so-called Roger's {\it Lie n-algebra of observables}   $L_{\infty}(M,\omega)$
(\cite{Rogers}). Referring to \cite{Rogers,RWZ,RW,Calliesetal, FLZ} for a full coverage of the relevant apparatus, not needed to full extent here,
we just point out that the latter is a graded vector space $L$ whose degree
$i$ pieces read
$$
\Lambda_{\rm Ham}^{n-1}  (M),\, \,  i=0, \qquad \Lambda^{n-1-i}(M), \,\,  i=1,2,\dots, n-1
$$
together with suitable multilinear maps denoted collectively by $\ell$.
The suffix ``Ham" refers to the {\it Hamiltonian (n-1)-forms}, i.e. those forms $H$ such that
$$
\iota_X \omega + d H =   \alpha(X) + dH = 0
$$
for a vector field $X$ preserving $\omega$ (i.e. ${\mathcal L}_X \omega = 0$), called, in turn,
a {\it Hamiltonian vector field} pertaining to
$H$.\par
A form $\beta$ is said to be {\it strictly} (resp. {\it globally}, resp. {\it locally}) conserved by an $\omega$-preserving vector field $X$ 
if ${\mathcal L}_X \beta = 0 $ (resp. ${\mathcal L}_X \beta$ is exact, resp. closed). Cartan's formula immediately shows
that closed forms are globally conserved; indeed, for such a form
$$
{\mathcal L}_{X} \beta =  d  \iota_{X} \beta  +  \iota_{X} d \beta  = d \iota_{X} \beta.
$$
Recall, from \cite{RWZ}, that  a {\it homotopy co-momentum map} is an $L_{\infty}$-algebra morphism - stemming from what is called an infinitesimal
action of ${\mathfrak g}$ on $M$ (with ${\mathfrak g}$ being the Lie algebra of a generic Lie group $G$, acting on $M$ by $\omega$-preserving vector fields)
$$
(f): {\mathfrak g} \to L_{\infty}(M, \omega)
$$
given explicitly by a sequence of linear maps
$$
(f)  = \{ f_i: \,\,\, \Lambda^i {\mathfrak g} \to \Lambda^{n-i}(M) \vert 0\leq i \leq n+1  \}
$$
fulfilling $f_0 = f_{n+1} = 0$ (we have tacitly set $\Lambda^{-1}(M) = 0$) and 
\begin{equation}\label{eq:f1_hcmm}
 \qquad \qquad {\rm Im} f_1 \in \Lambda^1_{\rm Ham}(M)
\end{equation}
together with (for $p \in \Lambda^k(\mathfrak{g})$):
\begin{equation}\label{eq:fk_hcmm}
-f_{k-1} (\partial p) = d f_k (p) + \varsigma(k) \iota(v_p) \omega
\end{equation}
($k=1,\dots n+1$). We explain the notation: first, if $ p = \xi_1 \wedge \xi_2 \wedge \dots \wedge \xi_k$, then $v_p = v_1 \wedge v_2 \wedge \dots \wedge v_k$  where $v_i \equiv v_{\xi_i}$ are the fundamental vector fields
associated to the action of $G$ on $M$. One sets $\iota(v_p) \omega = \iota(v_k)\dots\iota(v_1) \omega$, 
$\varsigma(k) := - (-1)^{\frac{k(k+1)}{2}}$ and defines $\partial \equiv \partial_k:  \Lambda^{k} {\mathfrak g} \to \Lambda^{k-1} {\mathfrak g}$  via
$$
\partial (\xi_1 \wedge \xi_2 \wedge \dots \wedge \xi_k) := \sum_{1\leq i< j \leq k} (-1)^{i+j}\, [\xi_i, \xi_j] \wedge \xi_1 \wedge \dots {\hat \xi}_i \wedge \dots \wedge {\hat \xi}_j \wedge \dots \xi_k
$$
(with $\,\hat{}\,$ denoting deletion as usual and with $\partial_0 = 0$; one has $\partial^2 = 0$). \par
Formula (\ref{eq:fk_hcmm}) tells us that the {\it closed} forms
$$
\mu_k := f_{k-1} (\partial p) +  \varsigma(k) \iota(v_p) \omega
$$
must actually be {\it exact}, with potential $-f_k(p)$.  Closure can be quickly ascertained as follows (in view of Lemma 2.16 in \cite{RWZ} and keeping in mind that $d\omega = 0$):
$$
\begin{matrix}
d (f_{k-1} (\partial p) +  \varsigma(k) \iota(v_p) \omega) & =  &  \varsigma(k) (-1)^k \iota(v_{\partial p})\omega -
\varsigma(k-1) \iota(v_{\partial p})\omega  \\
\phantom{void}& \equiv & [-\varsigma(k+1) -\varsigma(k-1)] \iota(v_{\partial p})\omega \\
\phantom{void}& = & 0 \\
\end{matrix}
$$
since in general $\varsigma(k)\varsigma(k+2) = -1$.  Notice that the special case $k=n+1$ asserts that the function $\mu_{n+1}(\cdot)$ is constant, and
its value is fixed by the condition 
$$
f_n(\partial p) + \varsigma(n+1)\iota(v_p)\omega = 0.
$$
 This can be rephrased, upon resorting to Section 9 in \cite{Calliesetal} (we use a different notation), by asserting that the following ${\mathfrak g}$-{\it (n+1)-cocycle} in the Chevalley-Eilenberg cochain (CE) complex $CE({\mathfrak g})$
$$
c_x  (p= \xi_1\wedge \xi_2 \wedge\dots \wedge \xi_{n+1}) = \iota_{v_p}\omega \!\mid_x
$$
ought to be a {\it boundary}: 
 $$
 c_x = \delta_{CE}(b)
 $$
 for a fixed but generic point $x \in M$ (the class $[c_x]$ being in general independent of $x \in M$, \cite{Calliesetal}, Cor. 9.3);  the operator $\delta_{CE}$ is the CE-differential defined by duality: $(\delta_{CE} \phi) (p) := \phi (\partial p)$, $\phi \in CE({\mathfrak g})$ and extended by linearity. Independence of $x$ is expressed via the formula (\cite{Calliesetal}, Prop. 9.1)
\begin{equation}\label{eq:c_chain}
 c_{x^{\prime}} - c_x =  \delta_{CE} (b)
\end{equation}
 where 
 $$ 
 b (v_1\wedge v_2 \wedge\dots \wedge v_{n}) := -\varsigma (n+1) \int_{\gamma} \iota (v_1\wedge v_2 \wedge\dots \wedge v_{n}) \omega
 $$
and $\gamma$ is a path connecting $x$ to $x^{\prime}$ (recall that $M$ is assumed to be connected).\par
We shall resume the above discussion in Subsection 2.3.\par

\subsection{The hydrodynamical Poisson bracket}
In the present Subsection we briefly review, for motivation and further applications, the symplectic geometrical portrait underlying the theory of perfect fluids, in its simplest instance.
We denote by ${\mathfrak g}$ the ({\it infinite dimensional}) Lie subalgebra of ${\mathfrak X}({\mathbb R}^3)$ consisting of the divergence-free vector fields on ${\mathbb R}^3$ 
that is the ``Lie algebra" of the ``Lie group"  $G = {\rm sDiff}\,({\mathbb R}^3)$ of volume preserving diffeomorphisms of ${\mathbb R}^3$.
As it is often done, we shall gloss over analytic subtleties
(see e.g. \cite{Arnold66,Arn-Khe,Eb-Mar,KM97} for more information).
We just recall here that $G$ is a {\it regular Lie group} in the sense of Kriegl-Michor (\cite{KM97}, 38.4) and that its associated exponential map is not even locally surjective (a quite general phenomenon).
We shall also tacitly assume that our fields {\it rapidly vanish at infinity}, so that convergence problems are avoided and boundary terms are absent.
The ``hydrodynamical" Lie bracket, equalling {\it minus} the standard one: $[\xi_1, \xi_2] = {\rm curl} (\xi_1 \times \xi_2)$ will be employed throughout.\par
Also, following e.g. \cite{Arn-Khe}, we shall consider the so-called  {\it regular dual} ${\mathfrak g}^*$ of ${\mathfrak g}$
consisting of all 1-forms modulo exact 1-forms:
$$
{\mathfrak g}^* := \Lambda^1({\mathbb R}^3)/ d \Lambda^0({\mathbb R}^3)
$$
together with the standard pairing ($\omega \in {\mathfrak g}^*$, $\xi \in {\mathfrak g} $)
$$
(\omega, \xi) = \int \langle \omega (x), \xi(x) \rangle \, d^3x .
$$
Nevertheless, we shall feel free to use suitable genuine distributional elements as well (i.e. {\it currents}, in the sense of de Rham, \cite{dR})
from the full topological dual (without introducing new notation for the latter). Everything will be clear from the context.\par
The (regular) dual ${\mathfrak g}^*$ is naturally interpreted as a Poisson manifold with respect to the hydrodynamical Poisson bracket  (Arnol'd --Marsden-Weinstein Lie-Poisson structure) ,
see e.g.\cite{Arn-Khe,Kuznetsov-Mikhailov80,MW83,Pe-Spe89,Pe-Spe92,Pe-Spe00,Spera16}:
$$
\{ F, G \} ([{\mathbf v}]) = \int_{ {\mathbb R}^3} \left\langle {\mathbf v}, \,  \left[\frac{\delta F}{\delta {\mathbf v}}, \frac{\delta G}{\delta {\mathbf v}}\right]  \right\rangle \,d^3x
$$
with ${\mathbf v} \in {\mathfrak g}$ (velocity field), ${\mathbf w} := {\rm curl}\,{\mathbf v}$,
its {\it vorticity},  with $[{\mathbf v}]$ denoting the ``gauge" class of ${\mathbf v}$: $[{\mathbf v}] = 
\{ {\mathbf v}+ \nabla f \}$.
The {\it Euler evolution}, reading, among others, in the so-called {\it vorticity form} 
$$
\frac{\partial {\mathbf w}}{\partial t} + [{\mathbf w}, {\mathbf v}] = 0
$$
is {\it volume preserving} and it also preserves the {\it symplectic leaves} of
${\mathfrak g}^*$ given by the $G$-{\it coadjoint orbits} ${\mathcal O}_{[{\mathbf v}]} \equiv {\mathcal O}_{{\mathbf w}}$.The symplectic structure on ${\mathcal O}_{\mathbf w}$ is the Kirillov-Kostant-Souriau (KKS)
(\cite{Kirillov01,Kostant70,Souriau70}):
$$
{\Omega_{KKS}}{([{\mathbf v}])} (ad^*_{{\mathbf b}}([{\mathbf v}]),ad^*_{{\mathbf c}})([{\mathbf v}]) = \int_{{\mathbb R}^3} \langle {\mathbf v}, [{\mathbf b},{\mathbf c}]\rangle \, d^3x=
\int_{{\mathbb R}^3} \langle {\mathbf w}, {\mathbf b} \times {\mathbf c}\rangle\, d^3x
$$
with the coadjoint action reading, explicitly, {\it up to a gradient} (not influencing calculations)
$$
ad^*_{{\mathbf b}}  ({\mathbf v}) = - {\mathbf w} \times {\mathbf b} 
\quad(\equiv ad^*_{{\mathbf b}}  ([{\mathbf v}]).
$$
The {\it Hamiltonian algebra} $\Lambda$ pertaining to ${\mathcal O}_{{\mathbf w}}$ consists of the so-called {\it Rasetti-Regge currents} originally introduced in \cite{Rasetti-Regge75} and further developed in
\cite{Pe-Spe89,Pe-Spe92,Pe-Spe00,Spera16,Bry}):
$$
\lambda_{\mathbf b} ({\mathbf v}) = \int\langle {\mathbf b},{\mathbf v} \rangle = \int\langle {\mathbf B},{\mathbf w} \rangle
$$
(with ${\rm curl}\, {\mathbf B} = {\mathbf b}$), fulfilling, for ${\mathbf b}$, ${\mathbf c} \in {\mathfrak g}$:
$$
\{\lambda_{\mathbf b}, \lambda_{\mathbf c} \} = \lambda_{[{\mathbf b},{\mathbf c}] },
$$
namely, the map
$$
{\mathfrak g} \ni {\mathbf b} \mapsto \lambda_{\mathbf b} \in \Lambda
$$
 is a {\it $G$-equivariant co-momentum map}
(observe in particular that
$\frac{\delta \lambda_{{\mathbf b}}}{\delta {\mathbf v}} = {\mathbf b}$).\par
\smallskip
The preceding portrait carries through to the {\it singular} vorticity case, in particular when the vorticity field is $\delta$-like and concentrated on a two-dimensional patch, filament or a loop.
Dealing with the latter case, we ultimately retrieve the {\it Brylinski manifold $Y$} consisting of {\it smooth oriented knots} (smooth embedded loops modulo orientation-preserving reparametrizations) together with its symplectic structure ${\Omega_Y}$ and the original Rasetti-Regge currents (see \cite{Bry,BeSpe06} for more details):
\begin{equation}
\Omega_Y (\cdot , \cdot) ({\gamma}) := \int_{\gamma} \nu (\dot{\gamma}, \cdot, \cdot) =
\int_{\gamma} \langle \dot{\gamma}, \cdot \times \cdot \rangle, \qquad \lambda_{\mathbf b} (\gamma) := \int_{\gamma} {\mathbf B}.
\end{equation}
Indeed, recall that, given a volume form $\nu$ on a 3-dimensional $M$, one gets, by {\it transgression}, a 2-form $\Omega$ on $LM$ via the formula
\begin{equation}
\Omega = \int_{S^1} ev^* (\nu )
\end{equation}
where $ev:  LM \times S^1\rightarrow M$ given by $ev (\gamma , t ) := \gamma(t)$ is the evaluation map of a loop $\gamma \in LM$
at a point $t \in S^1 \equiv [0,1] / \,\,{}_{\widetilde{}}\,$ (endpoint identification). More explicitly, given  tangent vectors $u$ and  $v$ at $\gamma$, the symplectic form reads
\begin{equation}
\Omega_{\gamma} (u , v)  = \int_{0}^1  \nu ({\dot \gamma}(t), u(t), v(t))\, dt
\end{equation}
(where we set ${\dot \gamma} = {d\gamma \over d t}$).
The above construction carries through to $Y$. In this case, the coadjoint orbits are labelled by the equivalence types of knots (via ambient isotopies),
by virtue of a result of Brylinski, see  \cite{Bry}.

\subsection{A hydrodynamical homotopy co-momentum map}
In this Subsection we elaborate on the previous discussion by introducing an explicit {\it homotopy co-momentum map} (HCMM), departing from the standard setting since our group $G = {\rm sDiff}\,({\mathbb R}^3)$ is infinite dimensional.
We start from the observation (\cite{CID}) that the volume form in ${\mathbb R}^3$, $\nu :=  dx \wedge dy \wedge dz $ can interpreted as a {\it multisymplectic form}: in this case
 the map $\alpha$  is bijective (in particular, injective). In coordinates, if $\xi = (\xi^i)$, then 
$$
\alpha(\xi) = \iota_{\xi} \nu = \xi^1 dy \wedge dz + \xi^2 dz \wedge dx + \xi^3 dx \wedge dy.
$$
Upon  introducing the Hodge $*$ relative to the standard Euclidean metric and the associated ``musical isomorphisms", we have ($\xi \in {\mathfrak{X}}({{\mathbb R}^3})$, $ \beta \in \Lambda^2({\mathbb R}^3)$):
$$
\alpha (\xi) =   * (\xi^{\flat}), \qquad \qquad \alpha^{-1}(\beta) =  (*\beta)^{\sharp}.
$$
Then we have, for $\xi \in \mathfrak {g}$ (via Cartan's formula)
$$
0 = {\mathcal L}_{\xi} \nu = d  \iota_{\xi} \nu  +  \iota_{\xi} d \nu  = d \iota_{\xi} \nu = {\rm div}(\xi) \nu
$$
and thus we have an isomorphism $\mathfrak {g} \cong Z^2({\mathbb R}^3)$ (closed 2-forms on ${\mathbb R}^3$). 
This will be important in the sequel. The above also expresses the fact that $\nu$ is a {\it strictly conserved} 3-form.\par

We shall now give the promised example of {\rm homotopy co-momentum map} emerging in fluid dynamics. Define, for $\mathbf{b} \in \mathfrak{g}$
$$
f_1(\mathbf{b}) := -  B
$$
where $B = \mathbf{B}^{\flat}$ and $ \mathbf{B}$ is again a vector potential for $\mathbf{b}$, i.e. ${\rm curl} \, \mathbf{B}= \mathbf{b}$, chosen e.g. in such a way that ${\rm div}\, \mathbf{B} = 0$ (Coulomb gauge).\par

It is immediately checked that
\begin{equation}\label{eq:condition1_hccm}
d f_1(\mathbf{b}) + \iota_{\mathbf{b}}\, \nu =    df_1 ({\mathbf b}) + \alpha({\mathbf b}) = 0.
\end{equation}
The above formula tells us that $ f_1(\mathbf{b})$ is a {\it Hamiltonian 1-form} for $\mathbf{b}$
(and, conversely, that the vector field $\mathbf{b}$ is a {\it Hamiltonian vector field} pertaining to
$ f_1(\mathbf{b})$, in accordance with (\ref{eq:f1_hcmm}) in Subsection 2.1.
Any $ f_1(\mathbf{b})$ above is also a  {\it Noether current} in the sense of Gotay et al., \cite{Gotay&al}.
In order to complete the definition of a homotopy co-momentum map, we just have  to find $f_2$,  satisfying formula (\ref{eq:fk_hcmm}) above.  
Indeed, for $k=1$ we retrieve (\ref{eq:condition1_hccm}).
The case $k=2$ reads
\begin{equation}\label{eq:condition2_hccm}
-f_1(\partial p) = df_2 (p) + \iota_{v_p} \nu.
\end{equation}
Let, for $\xi_i \in {\mathfrak g}$ ($i=1,2$),  $p = \xi_1 \wedge \xi_2$, so $\partial p = -[\xi_1, \xi_2]$ .

Then one checks that (using e.g. \cite{RWZ} Lemma 2.18, or the preceding subsection)
$$
df_1 ([\xi_1, \xi_2]) = d (\iota_{\xi_1 \wedge \xi_2} \nu) = -\iota_{[\xi_1,\xi_2]} \nu .
$$
(recall that $\iota_{\xi_1 \wedge \xi_2} \nu = \nu(\xi_1,\xi_2, \cdot)$).
Therefore, the 1-form $\mu_2(\xi_1,\xi_2) := f_1([\xi_1,\xi_2]) - \iota_{\xi_1 \wedge \xi_2} \nu$ is closed, hence exact, and (\ref{eq:condition2_hccm}) tells us that $f_2(p)$ is a potential for it and, as such, it is determined up to a constant $c(\xi_1,\xi_2)$.  In order to prove that we have a bona fide co-momentum map, we must have, in particular, for $q = \xi_1 \wedge \xi_2 \wedge \xi_3$, the explicit formula
\begin{equation}\label{eq:condition3_hccm}
	f_2(\partial q) = \nu(\xi_1, \xi_2, \xi_3) 
\end{equation}
which is a priori true up to a constant $c(\xi_1, \xi_2, \xi_3)$ by virtue of (\ref{eq:condition2_hccm}) and \cite{Calliesetal}, Lemma 9.2. However, the constant is in fact zero since 
 $\nu(\xi_1, \xi_2, \xi_3)$ vanishes at infinity and the same is true for $f_2(\partial q)$ upon solving the related Poisson equation
 $$
 \Delta f_2(\partial q) = \Delta \nu(\xi_1, \xi_2, \xi_3)
 $$
 (obtained via a straightforward computation; notice that we use the Riemannian Laplacian, which is {\it minus} the standard one). \par
An alternative derivation uses $x$-independence of the class $[c_x]$. Upon taking $S^3 = {\mathbb R}^3 \cup \{\infty \}$ , we have $c_{\infty } = 0$, hence
$c_x  = \delta_{CE} (b)$, with 
$$
b =  -\int_{\gamma_{\infty}} \iota (v_1\wedge v_2) \nu
$$
($\gamma_{\infty}$ being a path connecting $ x $ to ${\infty}$, cf. (2.3): the expression is meaningful in view of the assumed decay at infinity of our objects).
This is equivalent to the previous equation (\ref{eq:condition3_hccm}).
The function $f_2$ is the given by the solution of a Poisson equation
$$
f_2 =  \Delta^{-1} \delta \mu_2
$$
(in the present case $\delta = -* d * $). This completes the construction of the sought for homotopy co-momentum map.\par
\medskip
We define {\it Poisson brackets} (PB) via the expression
\begin{equation}
\{ f_1({\mathbf b}), f_1({\mathbf c}) \} (\cdot):= \iota_{\mathbf c} \iota_{\mathbf b} \nu (\cdot)= \nu({\mathbf b}, {\mathbf c}, \cdot)
\end{equation}
which will be employed below and in Section 4.\par
We may also naturally ask the question whether the above map $(f)$ is (infinitesimally) {\it G-equivariant}, in the sense of \cite{RWZ}: in particular, one should check the validity of the formula
$$
{\mathcal L}_{\xi} f_1({\mathbf b}) = f_1 ([\xi, {\mathbf b}])
$$
for all $\xi$, ${\mathbf b} \in {\mathfrak g}$. However, working out the two sides of the above equation easily yields, in particular, for $\xi = {\mathbf b}$, the equality
 $$
 d \, B({\mathbf b}) = 0,
 $$
that is, in vector terms $\langle {\mathbf B},{\mathbf b}\rangle = c = 0$ \ since ${\mathbf b}$ is compactly supported. However, if one considers a {\it flux tube} with
non zero {\it helicity} $\int \langle {\mathbf B},{\mathbf b} \rangle$ (see \cite{Moffatt-Ricca92,BeSpe06,Spe06} and below for further elucidation of this train of concepts), we get a contradiction. Notice that the argument does not depend on the choice of
$B$.
The lack of $G$-equivariance is not surprising, since our construction involves Riemannian geometric features.\par
We may now recap the preceding discussion via the following:
\begin{thm}
(i) The map $(f)$ previously given through the above $f_j: \Lambda^j {\mathfrak g} \to  \Lambda^{2-j} ({\mathbb R}^3)$, fulfilling (\ref{eq:condition1_hccm}),(2.5),(2.6),
yields a {\rm homotopy co-momentum map}; explicitly:
\begin{equation*}
f_1 = f_1(\mathbf{b}) := -  B, \quad (f_1 = \flat\circ{\rm curl}^{-1}); \qquad\quad f_2  
= \Delta^{-1} \delta \mu_2.
\end{equation*}
\par
(ii) The above HCCM transgresses, via the evaluation map ${\rm ev}: L{\mathbb R}^3 \times {\mathbb R} \ni (\gamma, t) \mapsto \gamma(t) \in {\mathbb R}^3$ to the hydrodynamical co-momentum map of Arnol'd and Marsden-Weinstein, defined on the Brylinski manifold $Y$ of oriented knots. \par
\smallskip
(iii)
Moreover, we have the formula
\begin{equation}\label{eq:bracketsformula}
\{ f_1({\mathbf b}), f_1({\mathbf c}) \}  - f_1([{\mathbf b},{\mathbf c}]) = -df_2 ({\mathbf b} \wedge {\mathbf c}).
\end{equation}

(iv) The map $(f)$ is {\rm not $G$-equivariant} in the sense of RWZ.
\end{thm}
{\bf Proof.} We have just to address point (ii): observe that the relevant piece of the homotopy co-momentum map is  $f_1$ which, under transgression, becomes
$$
\mu_{\mathbf b} = - \int_{\gamma} B = -\lambda_{\mathbf b}
$$
i.e., up to sign, the {\it Rasetti-Regge current} (RR) $\lambda_{\mathbf b}$ pertaining to ${\mathbf b} \in {\mathfrak g}$
independent of the choice of $B$.
This is in accordance with the general result in \cite{RWZ} asserting that, roughly speaking, homotopy co-momentum maps transgress to homotopy co-momentum maps on loop (and even mapping) spaces.
{\it Actually, the ansatz for $f_1$ term was precisely motivated by this phenomenon}.\par
Formula  (\ref{eq:bracketsformula}) in (iii) is just a rewriting of (\ref{eq:condition2_hccm}). \qed\par
\smallskip
\noindent
\subsection{A generalization to Riemannian manifolds}
 We ought to notice that a hydrodynamically flavoured homotopy co-momentum map can be similarly construed also for
an $(n+1)$-dimensional connected, compact, orientable Riemannian manifold $(M,g)$, upon taking its Riemannian volume form $\nu$ as a multisymplectic form and again the group $G$ of volume preserving diffeomorphism group as symmetry group. 
The divergence of a vector field $X$ is defined via ${\rm div}\, X := *d\!*\!X^{\flat} = -\delta X^{\flat} $; the operator $*$ is the Riemannian Hodge star (e.g. \cite{Warner,KM97}).
We can indeed prove the following result:
\begin{thm}
Let $(M,g)$ be a connected compact oriented Riemannian manifold of dimension $n+1$, $n\geq1$, with multisymplectic
form $\nu$ given by its Riemannian volume form, and such that the {\it de Rham}
cohomology groups $H_{dR}^{k}(M)$ vanish for $k=1,2,\dots n-1$ (one has necessarily $H_{dR}^{0} (M) = H_{dR}^{n+1} (M) = {\mathbb R}$). 
Let ${\mathfrak g}_0$ the Lie subalgebra of ${\mathfrak g}$ consisting of divergence-free vector fields vanishing at a point $x_0 \in M$. Then there exists an associated family of 
${\mathfrak g}_0$-homotopy co-momentum maps.\par 
\end{thm}
{\bf Proof.} 
As we have already noticed in general, the defining formula triggers a recursive construction starting from $f_1$, up to topological obstructions (we have a sequence of closed forms, which must be actually exact, together with the constraint $ f_n(\partial q) = (-1)^{\frac{(n+1)(n+2)}{2}}\nu(\xi_1,\dots\xi_{n+1})$, with $q= \xi_1 \wedge \dots \xi_{n+1}$, for the {\it constant} function $\mu_{n+1}(\cdot)$). In the present case, a natural candidate for the (n-1)-form $f_1$ can be readily manufactured via Hodge theory (see e.g. \cite{Warner}):
\begin{equation}\label{eq:f1_hcmm_riemann}
f_1(\xi) := -\Delta^{-1} \delta (\iota_{\xi} \nu)
\end{equation}
(the direct generalization of the preceding case)
after imposing $\delta f_1({\xi}) = 0$ (the analogue of the Coulomb gauge condition), provided one can safely invert the
Hodge Laplacian $\Delta = d\delta + \delta d$, this being the case if $H_{dR}^{n-1}(M) = 0$. One can of course alter the above definition by addition of an exact form. The 
topological assumptions made ensure that the entire procedure goes through unimpeded due to the formula
$$
df_k (\xi_1 \wedge\dots \wedge \xi_k) = \mu_k (\xi_1 \wedge\dots \wedge \xi_k), \qquad \qquad k=2,3,\dots n.
$$
One has finally to check that
$$
f_n (\partial (\xi_1 \wedge\dots \wedge \xi_{n+1})) =  -\varsigma(n+1) \iota (\xi_1 \wedge\dots \wedge \xi_{n+1})\nu
$$
but this is true once we observe that, since $c_{x_0} = 0$, the class $[c_x] = 0$ (cf. Subsection 2.1).
Therefore, we have, eventually, the compact formulae
$$
f_1(\xi) := -\Delta^{-1} \delta (\iota_{\xi} \nu); \quad\quad f_k = \Delta^{-1} \delta \mu_k, \qquad k=2\dots n.
$$

\hfill \qed
\par
\medskip 
{\bf Remark.} We notice that the above result holds, in particular, for  {\it homology spheres} such as, for instance, the celebrated Poincar\'e dodecahedral space. We point out that the case in which the intermediate homology groups are at most torsion (hence not detectable by de Rham techniques) is also encompassed: this is e.g. the case of {\it lens spaces}. Notice that $G$-equivariance cannot be expected a priori. Also notice that one could restrict to the natural symmetry group provided by the isometries of $(M,g)$. See e.g. \cite{RW15} for a general discussion of topological constraints to existence and uniqueness of homotopy co-momentum maps.\par

\subsection{Covariant phase space aspects}
We are now going to propose a multisymplectic interpretation of the hydrodynamical bracket which  
ties neatly with the topics discussed in previous sections, via {\it covariant phase space} \ ideas (\cite{KS,Gotay&al, FoRom,Zuckerman87,Crnkovic}, but without literally following the standard recipe, as we shall see shortly.
\par
Start with a 4-dimensional space-time $M =  {\mathbb R}^3 \times  {\mathbb R} \rightsquigarrow (x,y,z,t)$, define the obvious
trivial bundle
\begin{equation*}
 E = M \times \mathbb{R}^3  \to M
\end{equation*}
and interpret 
$\Sigma :={\mathbb R}^3  \rightsquigarrow (x,y,z)$ as a Cauchy ``submanifold" of $M$.\par
Any divergence-free vector field  can be viewed as an initial condition ${\mathbf v} (x,0)$ for the (volume-preserving) Euler evolution (at least for small times, but as we previously said, we do not insist on refined analytical nuances) ${\mathbf v} (x,t)$, yielding a section of $E$. Using the 3-volume form $\nu$, orienting fibres (notice that, when viewed on $E$, it is only pre-2-plectic, namely closed but degenerate), and
observing that we can set 
$$
{\mathcal J}^1 {\mathbf v}  := {\mathbf w} \quad (:= {\rm curl}\, {\mathbf v} )
$$
(the natural "covariant" jetification of the section ${\mathbf v} $-to be distinguished from the standard jetification $j^1$-if we wish to look at ${\mathbf v} $ as the vector space counterpart of a connection 1-form, yielding a section of the jet bundle $J^1E \to E$) we can rewrite the hydrodynamical bracket, mimicking 
\cite{FoRom}, as 
$$
\{ F, G \} ([{\mathbf v}]) = \int_{\Sigma = {\mathbb R}^3} \left\langle {\mathbf w}, \,  \frac{\delta F}{\delta {\mathbf v}} \times \frac{\delta G}{\delta {\mathbf v}}\right\rangle \,d^3x = \int_{\Sigma = {\mathbb R}^3} \nu ({\mathcal J}^1{\mathbf v}, \, \frac{\delta F}{\delta {\mathbf v}}, \frac{\delta G}{\delta {\mathbf v}}) \,d^3x =: (\star)
$$
since the variations $\frac{\delta F}{\delta {\mathbf v}}$ and $\frac{\delta G}{\delta {\mathbf v}}$ are {\it vertical} and divergence-free: $\delta F/ \delta {\mathbf v} = {\rm curl} \, (\delta F/ \delta {\mathbf w})$.
 Taking again ${\mathbf b} = {\rm curl}\, {\mathbf B}$ et cetera  and setting finally $F = \lambda_{\bullet}$ (see e.g. in particular \cite{Pe-Spe92,Spera16}) we see that the expression $(\star)$ can be manipulated to yield the expressive layout (with slight abuse of language)
$$
(\star) = \int_{\Sigma} ({{\mathcal J}^1}^*\nu) ({\mathbf v},{\mathbf B},{\mathbf C}) \,d^3x =  \int_{\Sigma}\nu ( {\mathcal J}^1{\mathbf v}, {\mathcal J}^1{\mathbf B}, {\mathcal J}^1{\mathbf C}) \,d^3x = \int_{\Sigma}\nu ( {\mathbf w}, {\mathbf b}, {\mathbf c}) \,d^3x
$$
(in full adherence with the discussion carried out in Section 2.2.).
The same portrait can be depicted, {\it mutatis mutandis}, for the singular case.
Ultimately, we reached the following conclusion: 
\begin{thm}
(i) The Poisson manifold ${\mathfrak g}^*$ can be naturally be interpreted as a (generalized) covariant phase space pertaining to the volume preserving Euler evolution: the latter indeed preserves the symplectic leaves of
${\mathfrak g}^*$ given by the $G$-coadjoint orbits ${\mathcal O}_{[{\mathbf v}]}$.\par
(ii) The above construction reproduces the symplectic structure of ${Y}$ upon taking singular vorticities, concentrated
on a smooth oriented knot: 
the covariant phase space picture is fully retrieved upon passing to
a 2-dimensional space-time $S^1 \times {\mathbb R}\rightsquigarrow (\lambda, t) $, with $\lambda \in S^1 \equiv \Sigma$ being a knot parameter (and staying of course with the same $\nu$). \par
\end{thm}
{\bf Remarks.} 1.  We stress the fact that we did not literally follow the standard "multisymplectic to covariant" recipe developed in \cite{FoRom}. In fact the multisymplectic manifold we consider is not the one prescribed by \cite{Gotay&al} since we directly took the standard volume form 
$\nu$ on ${\mathbb R}^3$ as a 2-plectic structure (or pre-2-plectic when pulled back to $E$), cf. \cite{CID}.
This neatly matches  Brylinski's theory  and fits with the stance long advocated, among others, by Rasetti and Regge and Goldin (see e.g. \cite{Rasetti-Regge75,Goldin71,Goldin87,Goldin12}, and \cite{Spera16} as well) pinpointing the special and ubiquitous role played by the group $G$.
Another motivation for considering $\nu$ is its pivotal role in the formulation of conservation theorems (see \cite{RWZ}). We shall pursue this aspect in what follows.\par
\smallskip
\noindent
2. In line with the preceding remark, notice that the above portrait can, in principle, be generalized to any volume form (on an orientable manifold), with its attached group $G$. The covariant phase space picture should basically persist in the sense that one might construct, in greater generality, an n-plectic structure out of an (n+1)-plectic one via an expression akin to $(\star)$. The (non) $G$-equivariance issue should be relevant in this context. \par
\section{A Hamiltonian 1-form for links}
We may specialize the considerations in Subsection 2.3 to the  case of links. 
The basic and quite natural idea is {\it to associate to a knot (or link ) a perfect fluid whose vorticity in concentrated thereon}
(cf. the preceding discussion on the Brylinski manifold).
As general references for knot theory we quote, among others, \cite{Rolfsen},
together with \cite{Bott-Tu82} for the algebraic-topological tools employed here. 
First recall  that in general the Poincar\'e dual $[\eta_S] \in H^{n-k}(M)$ of a k-dimensional closed oriented submanifold $S$ of an n-dimensional manifold $M$ is characterized by the property:
$$
\int_M \omega \wedge \eta_S = \int_S i^* \omega
$$
for any {\it closed, compactly supported k-form} $\omega$ on $M$ ($i: S \hookrightarrow M$ being the inclusion map). We shall indifferently view Poincar\'e duals as genuine forms or currents in the sense of de Rham (\cite{dR}).\par\par
\begin{figure}[h!]
\centering
\includegraphics[width=\textwidth]{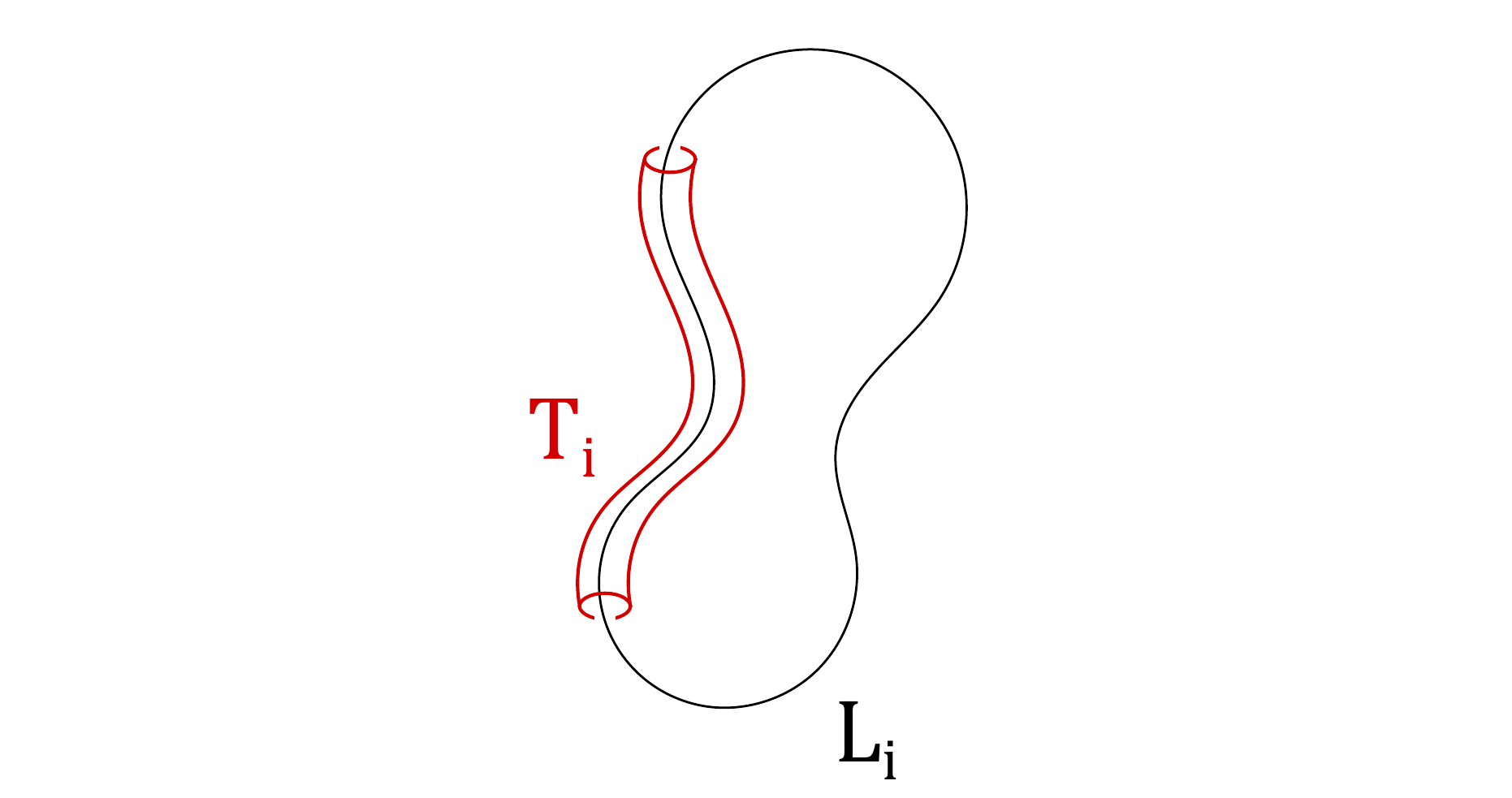}
\caption{Tubular neighbourhoods}
\label{fig: tubes}
\end{figure}

Building on 
\cite{Pe-Spe02,BeSpe06,Spe06}, let $ L = \cup_{i=1}^n L_i$ be an oriented link in ${\mathbb R}^3$ with components $L_i$, $i=1,\dots,n$ - required to be  {\it trivial} knots and let 
 $\omega_{ L_i}$ denote the Poincar\'e (or Thom) dual (class) associated to $L_i$: they are 2-forms localized in a  cross-section of a  suitable tubular neighbourhood $T_i$ around $L_i$ - with total fibre integral equal to one - see \cite{Bott-Tu82}, or, as currents, 2-forms which are $\delta$-like on $L_i$, see Figure \ref{fig: tubes}.\par

\begin{figure}[h!]
\centering
\includegraphics[width=\textwidth]{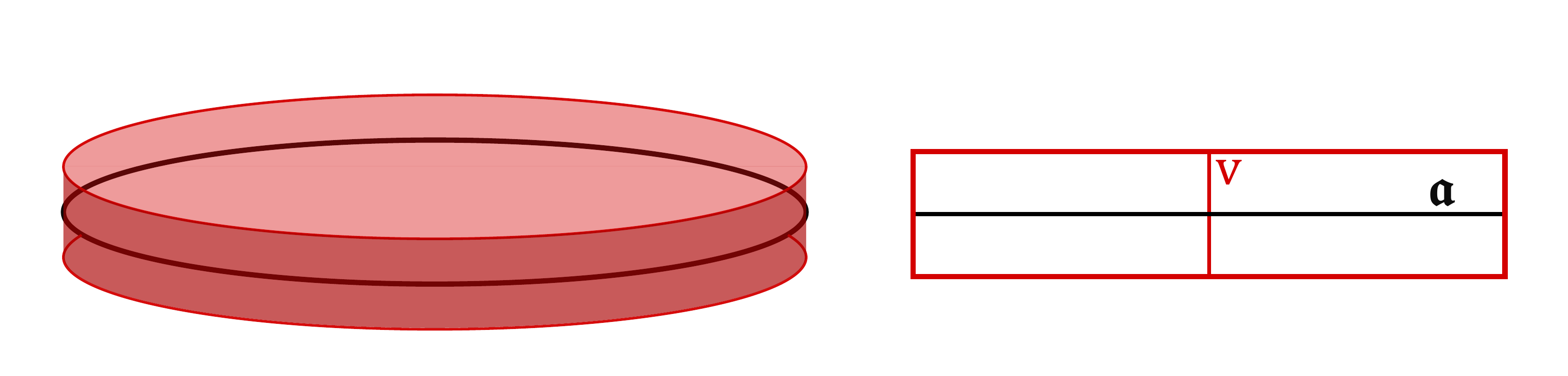}
\caption{Poincar\'e duals}
\label{fig: thom}
\end{figure}

Then take, for each $i=1,2,\dots,n$, a 1-form $v_{L_i}$ such that $dv_{L_i} = \omega_{L_i}$, namely, $v_{L_i} := \omega_{{\mathfrak a}_i}$ is the Poincar\'e dual (class) of a disc ${\mathfrak a}_i$ bounding
$L_i$ (a Seifert surface for the trivial knot $L_i$). Precisely:
$$
\partial {\mathfrak a}_i = L_i, \qquad \qquad dv_{L_i} = d\omega_{{\mathfrak a}_i} = \omega_{L_i} = \omega_{\partial {\mathfrak a}_i},
$$
see Figure \ref{fig: thom}.
We list, for the sake of clarity, (de Rham) cohomology and relative homology groups of $S^3 \setminus L$ with real coefficients, 
respectively, reading
$$
\begin{matrix}
H^0 (S^3 \setminus L) \cong H_3(S^3, L) \cong {\mathbb R} \\
\phantom{n} \\
\,\,H^1 (S^3 \setminus L) \cong H_2(S^3, L) \cong {\mathbb R}^{n} \\
\phantom{n} \\
\,\,\,\,\,\,\,H^2 (S^3 \setminus L) \cong H_1(S^3, L) \cong {\mathbb R}^{n-1} \\
\phantom{n} \\
\!\!\!H^3 (S^3 \setminus L) \cong H_0(S^3, L) \cong 0 \\
\end{matrix}
$$
The (de Rham classes of) the forms (or currents) $v_{L_i}$ generate in fact the cohomology group $H^1(S^3 \setminus L, {\mathbb R})$ (or, better, that of $S^3 \setminus T$, with $ T = \cup_{i=1}^n T_i$). Their homological counterparts are given by the (classes of) the discs ${\mathfrak a}_i$. One can also interpret the other groups:
in particular, elements in $H_1(S^3, L) $ can be represented by classes $[\gamma_{ij}]$ of (smooth) paths $\gamma_{ij}$ connecting two components $L_i$ and $L_j$, subject to the relation 
$[\gamma_{ij}] + [\gamma_{jk}] = [\gamma_{ik}]$.\par
Now set:
$$
\omega_{L} := \sum_{i=1}^n \omega_{L_i}
$$
(the vorticity 2-form for the link $L$) together with its velocity 1-form
$$
 v_{ L} = \sum_{i=1}^n v_{L_i}, \qquad \qquad  dv_{L} = \omega_{ L}.
$$
\begin{prop}
The position  
$$
H = v_{L}
$$
produces a {\it Hamiltonian 1-form} for links.
\end{prop}
 {\bf Proof.} The proof is straightforward: indeed for each component $L_i$, the Hamiltonian vector field $\xi_{L_i}$ for $v_{L_i} \equiv v_i$ is {\it minus} the vector field associated to the closed 2-form 
$\omega_{L_i}$ (via the map $\alpha$ of Section 2). Explicitly, one has (setting $\xi_L = \sum_{i=1}^n \xi_{L_i}$)
\begin{equation}\label{eq:poincaredual}
dv_L + \iota_{\xi_L} \nu = 0 
\end{equation}
\qed\par
\smallskip
{\bf Remark.}  Inspection of the very geometry of Poincar\'e duality shows that the velocity 1-forms $v_i$ correspond (upon approximation of the associated Euler equation) to the so-called LIA (Linear Induction Approximation) or  {\it binormal evolution}
of the ``vortex ring" $L_i$ (``orthogonal" to the discs ${\mathfrak a}_i$ - an easy depiction, cf. Figure \ref{fig: thom}), see \cite{Khe} for more information.
Formula (\ref{eq:poincaredual}) will be the prototype for the calculations in Section 4.\par

\smallskip
Let us define the {\it Chern-Simons (helicity) } 3-form:

$$
CS({L}) :=  v_{L} \wedge  \omega_{ L}.
$$
Integration of $CS({L})$ over ${\mathbb R}^3$ or $S^3$ yields an {\it integer} ${\mathcal H}(L) $, the {\it helicity} of $L$: 
$$
\int_{S^3} CS({L}) =: {\mathcal H}(L)  = \sum_{i,j=1}^n \ell(i,j),
$$
with $\ell(i,j) = \ell(j,i)$ being the {\it Gauss linking number} of components $L_i$ and $L_j$ if $i\neq j$ and $\ell(j,j)$ is the {\it framing} of $L_j$, equal to $\ell(L_j, L_j^{\prime})$ with $L_j^{\prime}$ being a section of the normal bundle of $L_j$, see e.g. \cite{Rolfsen,Moffatt-Ricca92,RN,Spe06} and below.
A regular projection of a link onto a plane produces a natural framing called the {\it blackboard framing}. \par %
On a Riemannian 3-manifold $M$, the helicity ${\mathcal H}$ pertaining to a perfect fluid with velocity
${\mathbf v}$ and vorticity ${\mathbf w} = {\rm curl}\,{\mathbf v}$ is given by (notation as in the previous sections)
$$
{\mathcal H} = \int_M \langle {\mathbf v},  {\mathbf w}\rangle = \int_M v \wedge dv
$$
(the last expression being the differential form counterpart).
Concretely, the helicity can be viewed as a measure of the mutual knotting of two generic flow lines, see \cite{Moffatt-Ricca92, Arn-Khe,Pe-Spe89,Pe-Spe92,Spe06} for a more extensive discussion. We used this concept in Subsection 2.3
to prove the non-equivariance of the hydrodynamical HCCM.\par
\section{A multisymplectic interpretation of Massey products }
In this section we resort to the techniques developed in Sections 2 and 3 above and propose a reformulation of the so-called
higher order linking numbers in multisymplectic terms. Ordinary and higher order linking numbers provide, among others,  a quite useful tool for the investigation of Brunnian phenomena in knot theory: recall that a link is {\it almost trivial} or {\it Brunnian} if upon removing any component therefrom one gets a trivial link.
They can be defined recursively in terms of Massey products, or equivalently, Milnor invariants,
by the celebrated Turaev-Porter theorem (see \cite{Fenn,Pe-Spe02,Spe06,Hebda-Tsau12}). We are going to review, briefly and quite concretely, the basic steps of the Massey procedure, read differential geometrically as in \cite{Pe-Spe02, Spe06,Hebda-Tsau12}, presenting at the same time our novel multisymplectic interpretation thereof. \par 
Let $L$ be an oriented link
with three or more components $L_j$. The cohomological reinterpretation of the ordinary linking number $\ell(1,2)$ of two components $L_1$ and $L_2$, say, starts from consideration of the closed 2-form
$$
\Omega_{12} := v_{1}\wedge v_{2}
$$
yielding the (integral) de Rham class 
$$
\langle L_1, L_2 \rangle :=   [\Omega_{12} ] \in H^2(S^3 \setminus L).
$$
The linking number $\ell(1,2)$ is non zero precisely when $\langle L_1, L_2 \rangle $, which, in $H_1(S^3,L)$ equals
$\ell(1,2) [\gamma_{12}]$, is non trivial.
If the latter class vanishes (i.e. $\Omega_{12}$ is exact), we have
\begin{equation}\label{eq:massey1}
dv_{12} + v_{1}\wedge v_{2}  = dv_{12} + \Omega_{12} = 0.
\end{equation}
for some 1-form $v_{12}$. Now, assuming that all the ordinary mutual linking numbers of the components under consideration vanish, one can manufacture
the (closed, direct check) 2-form (Massey product)
$$
\Omega_{123} =  v_{1}\wedge v_{23}  + v_{12}\wedge v_{3} 
$$
yielding a {\it third order linking number} (as a class):
$$
\langle L_1, L_2, L_3 \rangle :=   [\Omega_{123} ] \in H^2(S^3 \setminus L).
$$
If the latter class vanishes, we find a 1-form $v_{123}$ such that

\begin{equation}\label{eq:massey2}
dv_{123} + v_{1}\wedge v_{23}  + v_{12}\wedge v_{3}  =  dv_{123} + \Omega_{123}  = 0.
\end{equation}
It is then easy to devise a general pattern, giving rise to forms $v_I$, $\Omega_I$ ($I$ being a general multiindex).
Actually, everything can be organised - via Chen's calculus of iterated path integrals \cite{Chen,Chen4}- in terms of sequences of {\it nilpotent connections} ${\mathbf v}^{(k)}$, $k=1,2...$ on a trivial vector bundle over $S^3 \setminus L$ and their attached 
{\it curvature forms} ${\mathbf w}^{(k)}$ (ultimately, the $\Omega_I$, \cite{Pe-Spe02,Spe06,Tavares,Hain}), everything stemming from the {\it Cartan structure equation}
$$
d {\mathbf v}^{(k)} + {\mathbf v}^{(k)} \wedge {\mathbf v}^{(k)} = {\mathbf w}^{(k)}
$$
together with the ensuing Bianchi identity
$$
d {\mathbf w}^{(k)} +  {\mathbf v}^{(k)} \wedge {\mathbf w}^{(k)} - {\mathbf w}^{(k)} \wedge {\mathbf v}^{(k)} = 0
$$
(the latter implying closure of the forms $\Omega_I$).
In order to give a flavour of the general argument, start from the nilpotent connection ${\mathbf v}^{(1)}$
with its corresponding curvature ${\mathbf w}^{(1)}$:
\begin{equation*}
{\mathbf v}^{(1)} =
\begin{pmatrix}
0 & v_1 &  0 & 0 \\
0 & 0    & v_2   & 0    \\
0 & 0     &  0    &   v_3 \\
0 & 0     &  0    &   0 \\
\end{pmatrix}, \qquad
{\mathbf w}^{(1)} =
\begin{pmatrix}
0 & 0 &  \Omega_{12} = v_1\wedge v_2 & 0 \\
0 & 0    & 0   & \Omega_{23} = v_2\wedge v_3    \\
0 & 0     &  0    &   0 \\
0 & 0     &  0    &   0 \\
\end{pmatrix}
\end{equation*}
Then proceed similarly with
$$
{\mathbf v}^{(2)} =
 \begin{pmatrix}
0 & v_1 &  v_{12} & 0 \\
0 & 0    & v_2   & v_{23}    \\
0 & 0     &  0    &   v_3 \\
0 & 0     &  0    &   0 \\
\end{pmatrix},\qquad  {\mathbf w}^{(2)} =
\begin{pmatrix}
0 & 0 &  0 & \Omega_{123} = v_{1} \wedge v_{23} + v_{12} \wedge v_{3} \\
0 & 0    & 0   & 0    \\
0 & 0     &  0    &   0 \\
0 & 0     &  0    &   0 \\
\end{pmatrix}
$$
(we made use of $dv_{12} + \Omega_{12} = dv_{23} + \Omega_{23} = 0$), and so on.\par

\begin{figure}[h!]
\centering
\includegraphics[width=\textwidth]{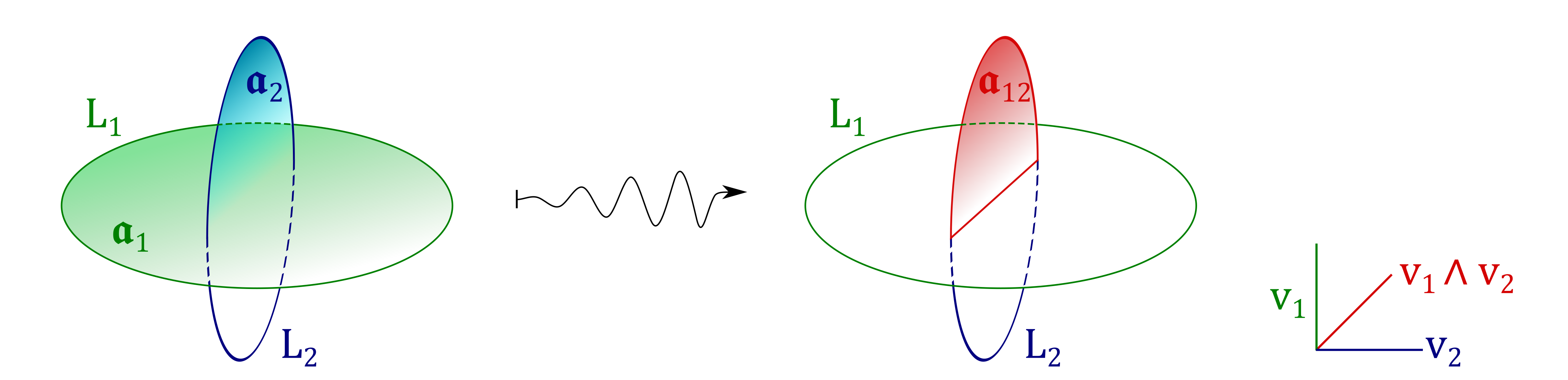}
\caption{Starting the Chen procedure}
\label{fig: chen} 
\end{figure}

\smallskip
Also recall that all forms $\Omega_{I}$ 
can be neatly interpreted, via Poincar\'e duality, as auxiliary (trivial) knots $L_I$, and $v_I$ as discs bounded by  $L_I$, in adherence to the considerations in Section 3, see \cite{Pe-Spe02,Spe06} for more details and worked out examples,
including the {\it Whitehead link} (involving fourth order linking numbers - with repeated indices) and the {\it Borromean rings} (exhibiting a third order linking number). Just notice here that, for instance, formula (\ref{eq:massey1}) becomes, intersection theoretically
$$
\partial {\mathfrak a}_{12} + {\mathfrak a}_1 \cap {\mathfrak a}_2 = 0,
$$
see Figure \ref{fig: chen}.
Formula (\ref{eq:massey2}) can be rewritten as
$$
dv_{123} + \iota_{\xi_{123}} \nu = 0 
$$
where ${\xi \equiv \xi_{123}}  = \alpha^{-1}(\Omega_{123})$. 
The above (``vorticity") vector field $\xi_{123}$ can be thought of as being concentrated on the knot corresponding to  $\xi_{123}$, or, alternatively, in a thin tube around it, when considering a bona fide Poincar\'e dual, cf. 
(\ref{eq:poincaredual}). \par
This tells us that $v_{123}$ is a {\it Hamiltonian 1-form} in the sense of \cite{RWZ} and the formula
$$
{\mathcal L}_{\xi} \Omega_{123} =  d  \iota_{\xi} \Omega_{123}  +  \iota_{\xi} d \Omega_{123}  = d \iota_{\xi} \Omega_{123} 
$$
 expresses the fact that $\Omega_{123}$ is a globally conserved 2-form, and the same holds for $\Omega_{12}$ and, in general, for $\Omega_I$,
with their corresponding vector fields $\xi_I$. Specifically, we have the following:\par

\begin{prop}
(i) The  volume form $\nu$ and all Massey 2-forms are globally conserved. \par
(ii) The 1-forms $v_{I} = f_1(\xi_I)$ are Hamiltonian with respect to the volume form.\par
\end{prop}

{\bf Proof.} Ad (i). This is clear since the mentioned forms are closed.\par
Ad (ii). The previous discussion can be carried out verbatim for a general  multiindex $I$:\par
$$
d v_I + \iota_{\xi_I} \nu = 0
$$
(an extension of (\ref{eq:poincaredual})), this yielding the second conclusion. \qed \par
The following is the main result of this section.
\begin{thm} With the above notation:\par
The 1-forms $v_I$ are {\rm first integrals in involution} with respect to the flow generated by the Hamiltonian vector field $\xi_{ L}$, namely
$$
{\mathcal L}_{\xi_{ L}} v_I = 0
$$
(i.e. the $v_I$'s are  {\rm strictly conserved})
and 
$$
\{v_I, v_J \} = 0
$$ 
(for multiindices $I$ and $J$).
\end{thm}
{\bf Proof.} Using Cartan's formula, we get
$$
{\mathcal L}_{\xi_{ L}} v_I = d \iota_{\xi_{ L}} v_I + \iota_{\xi_{ L}} d v_I = 
d \iota_{\xi_{ L}} v_I  -  \iota_{\xi_{ L}} \iota_{\xi_{ I}} \nu,
$$
but the second summand vanishes in view of the general expression
$$
\{v_\xi, v_\eta \}(\cdot) = \nu (\xi, \eta, \cdot)
$$
and of the peculiar structure of the vector fields involved (they either partially coincide or have disjoint supports). By the same argument,
one gets $\iota_{\xi_{ L}} v_I = 0$, in view of the Poincar\'e dual interpretation of $v_I$ (cf. Section 3), together with the second assertion;
a crucial point to notice  is that the auxiliary links obtained via Chen's procedure may be suitably split from their ascendants, this leading to
$$
\iota_{\xi_{ L}} v_I = 0,
$$
the consequent {\it strict} conservation of the $v_I$'s being then immediate.\par
\smallskip
Notice that, in particular, from
$$
\iota_{\xi_{ L}} v_L = 0
$$
(Poincar\'e dual interpretation again) we also get
$$
{\mathcal L}_{\xi_{ L}} v_L = 0
$$
(this is {\it not} to be expected a priori in multisymplectic geometry, cf.  \cite{RWZ}). \qed \par
\smallskip
\par
\smallskip
We ought to remark that, upon altering the $v_I$'s by an exact form, we may lose strict conservation, but in any
case global conservation is assured (the PB is an exact form, by (\ref{eq:bracketsformula}) in Section 2 and in view of commutativity of the
vector fields $\xi_I$ and $\xi_J$).\par
\smallskip
Ultimately, we can draw the conclusion that  {\it the Massey invariant route to ascertain the Brunnian character of a link can be mechanically understood as a recursive test of a kind of {\rm knot theoretic integrability}: the Massey linking numbers provide obstructions to the latter}. \par
\smallskip
Thus, somewhat curiously, higher order linking phenomena receive an interpretation in terms of  multisymplectic geometry, which is a sort of higher order symplectic geometry. Also, integrability comes in with a twofold meaning: first, higher order linking numbers emerge from the construction of a sequence of  flat, i.e. integrable nilpotent connections; second, this very process yields first integrals in involution in a mechanical sense. \par
\section{Conclusions and outlook}
In this note we constructed a homotopy co-momentum map in a hydrodynamical context, whereby we gave, as an application, a multisymplectic reinterpretation of the Massey higher order linking numbers, together with an extension thereof in a Riemannian geometric framework. 
We have also exhibited a covariant phase space interpretation of the geometrical framework involved. \par
The multisymplectic approach appears to be very promising for further advancement in this area. Also, the notion of integrability cropping up in our analysis of Massey products may deserve further scrutiny in a general multisymplectic context. 
We hope to tackle (at least some of) the open questions raised in this paper elsewhere.
\par
\medskip
\noindent
{\bf Acknowledgements.}
The authors, both members of the GNSAGA group of INDAM, acknowledge support from Unicatt local D1-funds (ex MIUR 60\% funds). They are indebted to T. Wurzbacher and M. Zambon for enlightening discussions. They are also grateful to Marcello Spera for help with graphics. \par

\end{document}